\documentclass[final]{siamltex}

\usepackage{color}
\usepackage{amsmath,amssymb,mathrsfs,eucal}
\usepackage[dvips,colorlinks,hypertex]{hyperref}

\newtheorem{rmk}[theorem]{Remark}
\newenvironment{remark}{\begin{rmk}\rm\!\!}{\end{rmk}}

\title{When do nonlinear filters achieve\\ maximal accuracy?\thanks{The 
author is partially supported by the NSF RTG Grant DMS-0739195.}}

\author{Ramon van Handel\thanks{Department of Operations Research and 
Financial Engineering, Princeton University, Princeton, NJ 08544
({\tt rvan@princeton.edu}).}}

\begin{document}

\maketitle

\begin{abstract}
The nonlinear filter for an ergodic signal observed in white noise is said 
to achieve maximal accuracy if the stationary filtering error vanishes as 
the signal to noise ratio diverges.  We give a general characterization of 
the maximal accuracy property in terms of various systems theoretic 
notions.  When the signal state space is a finite set explicit necessary 
and sufficient conditions are obtained, while the linear Gaussian case 
reduces to a classic result of Kwakernaak and Sivan (1972).
\end{abstract}

\begin{keywords} 
nonlinear filtering, maximal accuracy, systems theory, small noise limit
\end{keywords}

\begin{AMS}
93E11, 60G10, 62M20, 93B07, 94A12
\end{AMS}

\pagestyle{myheadings}
\thispagestyle{plain}
\markboth{R. VAN HANDEL}{WHEN DO NONLINEAR FILTERS ACHIEVE MAXIMAL 
ACCURACY?}


\section{Introduction}

Let $(X_t)_{t\ge 0}$ be a signal of interest, which we model as an
ergodic Markov process.  It is often the case that the detection of such a 
signal is imperfect: only some function of the signal may be directly 
observable, and the observations are additionally corrupted by additive 
white noise.  That is, one observes in practice the integrated observation 
process
$$
	Y_t = \int_0^t h(X_s)\,ds + \kappa\,B_t,
$$
where $B_t$ is a Wiener process independent of $(X_t)_{t\ge 0}$, 
$\kappa>0$ determines the strength of the corrupting noise, and $h$ is a 
(possibly nonlinear and noninvertible) function of the signal.  When only 
the imperfect observations $(Y_s)_{s\le t}$ are available, the exact 
value of the signal $X_t$ can certainly not be detected with 
arbitrary precision, even when $t$ is very large (so that we have a long 
observation history at our disposal).

To improve the accuracy of our detector, we must decrease the strength of 
the corrupting noise.  It is intuitively obvious that as $\kappa\to 0$, we 
will eventually be able to determine precisely the value of $h(X_t)$.  
However, when the function $h$ is not invertible (as is the case in many 
engineering systems of practical interest), this does not necessarily 
imply that we will be able to determine precisely the value of the signal 
itself.  The optimal estimate of the signal $X_t$, given the observation 
history $(Y_s)_{s\le t}$, is called the nonlinear filter.  We say that the 
filter \emph{achieves maximal accuracy} if, as the noise strength 
vanishes, the stationary filtering error vanishes also---i.e., if as 
$t\to\infty$ and $\kappa\to 0$, we are able to determine precisely the 
value of the signal.

When do nonlinear filters achieve maximal accuracy?  In the special linear 
Gaussian case, where the nonlinear filter reduces to the Kalman-Bucy 
filter, this question was first posed and resolved in a well known paper 
of Kwakernaak and Sivan \cite{KS72}.  Somewhat surprisingly, the answer is 
far from trivial and the proof given by Kwakernaak and Sivan is reasonably 
involved.  In fact, Kwakernaak and Sivan chiefly study the dual 
deterministic control problem with `cheap' control. Their proof is not 
probabilistic in nature, but is based on a delicate analysis of the 
associated Riccati equation.

Very little appears to be known beyond the linear Gaussian case.  To the 
best of our knowledge the only nonlinear result is due to Zeitouni and 
Dembo \cite{ZD88}, who study a special class of diffusion signals with 
nonlinear drift term and linear observations.  Their result, however, also 
reduces to the linear Gaussian case: the key step in the proof is to 
estimate the filtering error by that of an auxiliary Kalman-Bucy filter.

The purpose of this paper is to investigate the maximal accuracy problem 
in a general setting.  After setting up the problem and introducing the 
relevant concepts in section \ref{sec:prelim}, we proceed in section 
\ref{sec:main} to relate the maximal accuracy property of the filter to 
several systems theoretic notions (theorem \ref{thm:main} below).  The 
proof of our main result follows from simple probabilistic arguments.  
Then, in section \ref{sec:finst}, we apply our general result to provide a 
complete characterization of the maximal accuracy property for the case 
where the signal is a finite state Markov process.  The resulting 
necessary and sufficient condition---observability of the model after time 
reversal, together with a condition of the graph coloring type---is easily 
verified, but is surprisingly quite different in nature than the result 
for linear Gaussian systems.

Finally, in section \ref{sec:ks}, we revisit the linear Gaussian setting 
and provide a complete proof of the result of Kwakernaak and Sivan using 
our general characterization.  Though this does not lead to new results, 
our approach does not use the explicit form of the filtering equations and 
some parts of the proof are significantly simpler than that of 
\cite{KS72}.  We believe that our approach takes a little of the mystery 
out of the result of Kwakernaak and Sivan by placing it within a general 
probabilistic framework.

\subsection*{Acknowledgment}

The problem studied in this paper was posed to me by Prof.\ Ofer Zeitouni 
during a visit to the University of Minnesota in October 2008.  I am 
indebted to him for arranging this visit and for our many subsequent 
discussions on this topic, without which this paper would not have been 
written.

\section{Preliminaries}
\label{sec:prelim}

We suppose that defined on a probability space 
$(\Omega,\mathcal{F},\mathbf{P})$ is a stationary Markov process 
$(X_t)_{t\in\mathbb{R}}$ with c{\`a}dl{\`a}g sample paths in the Polish 
state space $E$, and we denote its stationary measure as 
$\mathbf{P}(X_0\in A)=\pi(A)$.  Moreover, we presume that the 
probability space supports an $n$-dimensional two-sided Wiener process 
$(B_t)_{t\in\mathbb{R}}$ that is independent of 
$(X_t)_{t\in\mathbb{R}}$.  Let us define for every $\kappa>0$ the 
$\mathbb{R}^n$-valued observation process $(Y_t^\kappa)_{t\in\mathbb{R}}$ 
according to the expression
$$
	Y_t^\kappa = \int_0^t h(X_s)\,ds + \kappa\,B_t,
$$
where $h:E\to\mathbb{R}^n$ is a given observation function.  
$X_t$ is called the signal process, and $Y_t$ is called the observations 
process.  In addition, we introduce the following notation. Let $\tilde 
X_t=X_{-t}$ be the time reversed signal, and note that $\tilde X_t$ is 
again a stationary Markov process under $\mathbf{P}$ with invariant 
distribution $\pi$.  We denote
\begin{equation*}
\begin{array}{ll}
	\mathcal{F}_{I}^X=\sigma\{X_s:s\in I\},&\qquad\quad
	\mathcal{F}_{I}^{\tilde X}=\sigma\{\tilde X_s:s\in I\},\\
	\mathcal{F}_{I}^{h(X)}=\sigma\{h(X_s):s\in I\},&\qquad\quad
	\mathcal{F}_{I}^{h(\tilde X)}=\sigma\{h(\tilde X_s):s\in I\}
\end{array}
\end{equation*}
for $I\subset\mathbb{R}$, while
$$
	\mathcal{F}_{[a,b]}^{Y,\kappa}=\sigma\{Y_s^\kappa-Y_a^\kappa:s\in 
		[a,b]\},\qquad
	\mathcal{F}_{[a,\infty[}^{Y,\kappa}=
		\bigvee_{b>a}\mathcal{F}_{[a,b]}^{Y,\kappa},\qquad
	\mathcal{F}_{]-\infty,b]}^{Y,\kappa}=
		\bigvee_{a<b}\mathcal{F}_{[a,b]}^{Y,\kappa}
$$
for $a\le b$.  Finally, for any probability measure $\mu\ll\pi$, we define
$$
	\mathbf{P}^\mu(A) = \mathbf{E}\left(I_A~
	\frac{d\mu}{d\pi}(X_0)
	\right).
$$
Note that $\mathbf{P}^\mu(X_0\in A)=\mathbf{P}^\mu(\tilde X_0\in 
A)=\mu(A)$ by construction, while 
$$
	\mathbf{E}^\mu(f(X_t)|\mathcal{F}^X_{]-\infty,s]}) =
	\mathbf{E}(f(X_t)|\mathcal{F}^X_{]-\infty,s]}),
	\qquad
	\mathbf{E}^\mu(f(\tilde X_t)|\mathcal{F}^{\tilde X}_{]-\infty,s]}) =
	\mathbf{E}(f(\tilde X_t)|\mathcal{F}^{\tilde X}_{]-\infty,s]})
$$
for any $t\ge s\ge 0$ and bounded function $f$ by the Bayes formula.  
Therefore, under $\mathbf{P}^\mu$, the one-sided processes $(X_t)_{t\ge 
0}$ and $(\tilde X_t)_{t\ge 0}$ are still Markov with the same transition 
probabilities as under $\mathbf{P}$, but with the initial measure $\mu$ 
instead of $\pi$.

\begin{remark}
If we are given a transition semigroup for the Markov process $(X_t)_{t\ge 0}$,
then we can construct $\mathbf{P}^\mu|_{(X_t,Y_t)_{t\ge 0}}$ even for 
$\mu\not\ll\pi$.  However, the transition semigroup for the time reversed 
process $(\tilde X_t)_{t\ge 0}$ under $\mathbf{P}$ is defined implicitly 
only up to $\pi$-a.s.\ equivalence in terms of the regular conditional 
probabilities $\mathbf{P}(X_{-t}\in\,\cdot\,|X_0)$.  Therefore, for the 
time reversed process, we can not unambiguously define 
$\mathbf{P}^\mu(\tilde X_t\in A)$ for $\mu\not\ll\pi$.  We will therefore 
restrict our attention throughout to probability measures $\mu\ll\pi$, 
except in remark \ref{rem:ergodic} and lemmas \ref{lem:erg1} and 
\ref{lem:erg2} below where only $(X_t,Y_t)_{t\ge 0}$ under 
$\mathbf{P}^\mu$ is considered (and not the time reversed part).
\end{remark}

The \textit{nonlinear filter} of the signal $X_t$ given the noisy 
observations $Y_s^\kappa$, $0\le s\le t$ is defined as the regular 
conditional probability $\mathbf{P}(X_t\in\,\cdot\,|
\mathcal{F}^{Y,\kappa}_{[0,t]})$.  By construction, the filter minimizes 
the mean square estimation error $e_t(f,\kappa)$ for every test function 
$f:E\to\mathbb{R}$ with $\int f^2d\pi<\infty$: i.e.,
$$
	e_t(f,\kappa) = \mathbf{E}\left(\left\{
	f(X_t)-\mathbf{E}\big(f(X_t)\big|\mathcal{F}^{Y,\kappa}_{[0,t]}
	\big)
	\right\}^2\right)
$$
is minimal ($e_t(f,\kappa)\le \mathbf{E}(\{f(X_t)-\hat F\}^2)$ for every
$\mathcal{F}^{Y,\kappa}_{[0,t]}$-measurable $\hat F$).

\begin{lemma}
For every test function $f$ with $\int f^2d\pi<\infty$ and every noise 
strength $\kappa\ge 0$, the mean square error $e_t(f,\kappa)$ converges to 
the stationary error
$$
	\lim_{t\to\infty}e_t(f,\kappa) = 
	\mathbf{E}\left(\left\{
	f(X_0)-\mathbf{E}\big(f(X_0)\big|\mathcal{F}^{Y,\kappa}_{]-\infty,0]}
	\big)
	\right\}^2\right) := e(f,\kappa).
$$
\end{lemma}

\begin{proof}
The result follows directly from the stationarity of $\mathbf{P}$ and the
martingale convergence theorem.
\qquad
\end{proof}

Our interest is in the behavior of $e(f,\kappa)$ in the limit $\kappa\to 
0$ where the observation noise vanishes.  In particular, we aim to 
understand when the filter achieves \textit{maximal accuracy}.

\begin{definition}
\label{def:maxac}
The filter is said to achieve \emph{maximal accuracy} if 
$e(f,\kappa)\to 0$ as $\kappa\to 0$ whenever $\int f^2d\pi<\infty$, i.e., 
if the true location of the signal is revealed in the stationary limit 
when the observation noise vanishes asymptotically.
\end{definition}

Our main results relate the maximal accuracy problem to certain structural 
properties of a systems theoretic flavor.  We define these notions presently.

\begin{definition}
\label{def:recon}
The model is said to be \emph{reconstructible} if 
$$
	\mu,\nu\ll\pi\quad\mbox{and}\quad
	\mu\ne\nu
	\qquad\mbox{implies}\qquad
	\mathbf{P}^\mu|_{\mathcal{F}^{h(X)}_{]-\infty,0]}}
	\ne \mathbf{P}^\nu|_{\mathcal{F}^{h(X)}_{]-\infty,0]}}.
$$
\end{definition}

\begin{definition}
\label{def:srecon}
The model is said to be \emph{strongly reconstructible} if 
$$
	\mu,\nu\ll\pi\quad\mbox{and}\quad
	\mu\perp\nu
	\qquad\mbox{implies}\qquad
	\mathbf{P}^\mu|_{\mathcal{F}^{h(X)}_{]-\infty,0]}}
	\perp \mathbf{P}^\nu|_{\mathcal{F}^{h(X)}_{]-\infty,0]}}.
$$
\end{definition}

\begin{definition}
\label{def:inv}
The model is said to be \emph{invertible} if for every $t\ge s$, the
random variable $X_t$ coincides $\mathbf{P}$-a.s.\ with a 
$\sigma(X_s)\vee\mathcal{F}^{h(X)}_{[s,t]}$-measurable random variable.
\end{definition}

\begin{definition}
\label{def:sinv}
The model is said to be \emph{stably invertible} if 
for every $t$ the random variable $X_t$ coincides $\mathbf{P}$-a.s.\ with 
a $\mathcal{F}^{h(X)}_{]-\infty,t]}$-measurable random variable.
\end{definition}

Let us finish this section with some remarks on these definitions.

\begin{remark} 
In the setting of deterministic linear systems theory, the notion of 
reconstructibility dates back to Kalman \cite{KFA69}, see also 
\cite{KS72bk}.  Our definition \ref{def:recon} in the stochastic setting
is close to a similar notion that plays an important role in the 
realization theory of stationary Gaussian processes \cite{Ruc78,LP85}.  
Reconstructibility is essentially the time reversed counterpart of the 
notion of observability \cite{Van09}, though as discussed above we must 
restrict to probability measures $\mu,\nu\ll\pi$.
\end{remark}

\begin{remark}
\label{rem:stabinv}
By the stationarity of $\mathbf{P}$, definitions \ref{def:inv} and 
\ref{def:sinv} can be restricted without loss of generality to the case 
$t=0$.  Thus invertibility means $X_0$ can be written as a function of 
$X_s$ at some previous time $s$ and all the intermediate noiseless 
observations $h(X_r)$, i.e., $X_0 = F_s[X_s,(h(X_r))_{s\le r\le 0}]$ for 
any $s<0$.  This idea is well known in the deterministic setting; see, 
e.g., \cite{Moy77} in the linear case and \cite{Hir79} in the nonlinear 
case. We think of the inverse $F_s$ as being `stable' if it becomes 
independent of $X_s$ as $s\to -\infty$; it therefore makes sense to talk 
of stable inversion when $X_0$ can be written as a function 
$X_0=F_{-\infty}[(h(X_r))_{r\le 0}]$ of all past noiseless observations.
\end{remark}

\begin{remark}
Suppose that the model is invertible.  Then certainly $X_0$ is 
$\mathbf{P}$-a.s.\ $\sigma\{X_{u},h(X_r):u\le s,~r\le 0\}$-measurable for 
every $s<0$.  In particular,
$$
	X_0\quad\mbox{is}\quad
	\mathbf{P}\mbox{-a.s.}\quad
	\bigcap_{s\le 0}\Big(\mathcal{F}^{X}_{]-\infty,s]}\vee
	\mathcal{F}^{h(X)}_{]-\infty,0]}\Big)\mbox{-measurable}.
$$
Now suppose that the left tail $\sigma$-field $\bigcap_{s\le 0}
\mathcal{F}^{X}_{]-\infty,s]}$ is $\mathbf{P}$-trivial, i.e.,
the signal is ergodic in a weak sense \cite{Tot70}.  Then it is tempting 
to exchange the order of intersection and supremum, as follows:
$$
	X_0\mbox{ is }\mathbf{P}\mbox{-a.s.\ }
	\bigcap_{s\le 0}\Big(\mathcal{F}^{X}_{]-\infty,s]}\vee
	\mathcal{F}^{h(X)}_{]-\infty,0]}\Big)\stackrel{?}{=}
	\Bigg[\bigcap_{s\le 0}\mathcal{F}^{X}_{]-\infty,s]}\Bigg]\vee
	\mathcal{F}^{h(X)}_{]-\infty,0]}=
	\mathcal{F}^{h(X)}_{]-\infty,0]}\mbox{-measurable}.
$$
This would indicate that invertibility plus ergodicity implies stable 
invertibility.  However, the exchange of intersection and supremum is not 
necessarily permitted, as an illuminating counterexample in \cite{BCL04} 
shows.  This conclusion is therefore invalid.  A further discussion of 
this problem can be found in \cite{vW83}.  In particular, it is evident 
that the present problem is closely related to the \textit{innovations 
problem} which is discussed in \cite{vW83}.  Another closely related 
problem, that of the stability of the nonlinear filter, is discussed in 
detail in \cite{Van08}; however, it should be noted that the 
nondegeneracy assumption made there is manifestly absent in the problems 
discussed here.
\end{remark}

\begin{remark}
\label{rem:ergodic}
Suppose the signal is not started in the stationary distribution, but 
in a distribution $\mu$ that is not even necessarily absolutely continuous 
with respect to $\pi$.  In this setting, the maximal achievable accuracy 
problem is to determine whether
$$
	e_t^\mu(f,\kappa) = \mathbf{E}^\mu\left(\left\{
	f(X_t)-\mathbf{E}^\mu\big(f(X_t)\big|\mathcal{F}^{Y,\kappa}_{[0,t]}
	\big)
	\right\}^2\right)
$$
converges to zero as $t\to\infty$, $\kappa\to 0$ (in that order).  
However, we will presently show that if $\|\mathbf{P}^\mu(X_t\in\,\cdot\,)
-\pi\|_{\rm TV}\to 0$ as $t\to\infty$, this problem reduces to the 
stationary problem where $\mu=\pi$.  In particular, \emph{when the signal 
is ergodic, the maximal achievable accuracy problem always reduces to the 
stationary case} (by ergodic we mean $\|\mathbf{P}^\mu(X_t\in\,\cdot\,)
-\pi\|_{\rm TV}\to 0$ for all $\mu$).  This strongly motivates our choice 
to study directly the stationary problem in the remainder of this paper.

Let us now make these claims precise in the form of two lemmas.  For 
simplicity, we concentrate on bounded functions, which is not a 
significant restriction.
\end{remark}

\begin{lemma}
\label{lem:erg1}
Suppose that $f$ is a bounded measurable function.  Then for any $\kappa$
$$
	\|\mathbf{P}^\mu(X_t\in\,\cdot\,)-\pi\|_{\rm TV}
	\xrightarrow{t\to\infty} 0\qquad\mbox{implies}\qquad
	\limsup_{t\to\infty}\,e_t^\mu(f,\kappa)\le e(f,\kappa).
$$
Thus $e(f,\kappa)\to 0$ as $\kappa\to 0$ implies
$e_t^\mu(f,\kappa)\to 0$ as $t\to\infty$, $\kappa\to 0$ (in that order).
\end{lemma}

\begin{lemma}
\label{lem:erg2}
Suppose that $\|\mathbf{P}^\mu(X_t\in\,\cdot\,)-\pi\|_{\rm TV}\to 0$ for 
all $\mu$ (i.e., the signal is ergodic). Then $e_t^\mu(f,\kappa)\to 
e(f,\kappa)$ for all $\kappa>0$, $\mu$, and bounded measurable $f$.
\end{lemma}

{\em Proof of lemmas \ref{lem:erg1} and \ref{lem:erg2}}.
Let $P_t$ be the Markov semigroup of the signal $\mu P_t=
\mathbf{P}^\mu(X_t\in\,\cdot\,)$.  We basically follow Kunita 
\cite{Kun71}.  First, by Jensen's inequality
\begin{equation*}
\begin{split}
	e_{t+s}^\mu(f,\kappa) &= 
	\mathbf{E}^\mu(f(X_{t+s})^2)-
	\mathbf{E}^\mu\!\left(\mathbf{E}^\mu\bigg(
	\mathbf{E}^\mu\big(f(X_{t+s})\big|\mathcal{F}^{Y,\kappa}_{[0,t+s]}
	\big)^2\bigg|\mathcal{F}^{Y,\kappa}_{[s,t+s]}\bigg)\right) \\
	&\le 
	\mathbf{E}^\mu(f(X_{t+s})^2)-
	\mathbf{E}^\mu\!\left(
	\mathbf{E}^\mu\big(f(X_{t+s})\big|\mathcal{F}^{Y,\kappa}_{[s,t+s]}
	\big)^2\right) \\
	&=
	\mathbf{E}^{\mu P_s}(f(X_t)^2)-
	\mathbf{E}^{\mu P_s}\!\left(
	\mathbf{E}^{\mu P_s}\big(f(X_t)\big|\mathcal{F}^{Y,\kappa}_{[0,t]}
	\big)^2\right)
	= e_t^{\mu P_s}(f,\kappa)
\end{split}
\end{equation*}
for $0<s<t$.  Now suppose that $\|\mu P_s-\pi\|_{\rm TV}\to 0$. Then it 
follows trivially (as $f$ is bounded) that 
$\mathbf{E}^{\mu P_s}(f(X_t)^2)\to\mathbf{E}(f(X_t)^2)$.  On the other 
hand, we can estimate
\begin{multline*}
	\left|
	\mathbf{E}^{\mu P_s}\!\left(
	\mathbf{E}^{\mu P_s}\big(f(X_t)\big|\mathcal{F}^{Y,\kappa}_{[0,t]}
	\big)^2\right) 
	-
	\mathbf{E}\!\left(
	\mathbf{E}\big(f(X_t)\big|\mathcal{F}^{Y,\kappa}_{[0,t]}
	\big)^2\right)\right|
	\\ \mbox{}
	\le
	\|f\|_\infty^2\|\mu P_s-\pi\|_{\rm TV}
	+ 2\|f\|_\infty
	\mathbf{E}\!\left(
	\left|
	\mathbf{E}^{\mu P_s}\big(f(X_t)\big|\mathcal{F}^{Y,\kappa}_{[0,t]}
	\big)-
	\mathbf{E}\big(f(X_t)\big|\mathcal{F}^{Y,\kappa}_{[0,t]}
	\big)\right|\right).
\end{multline*}
Using the identity $\mathbf{E_P}(|\mathbf{E_P}(X|\mathcal{G})-
\mathbf{E_Q}(X|\mathcal{G})|)\le 4\|X\|_\infty
\|\mathbf{P}-\mathbf{Q}\|_{\rm TV}$ \cite[theorem 3.1]{CP05}, it follows 
that the left hand side of this expression converges to zero.  We have 
therefore shown that if $\|\mu P_s-\pi\|_{\rm TV}\to 0$, then
$e_t^{\mu P_s}(f,\kappa)\to e_t(f,\kappa)$ as $s\to\infty$.  In 
particular,
\begin{multline*}
	\limsup_{t\to\infty} e_t^\mu(f,\kappa) =
	\limsup_{t\to\infty}\limsup_{s\to\infty}e_{t+s}^\mu(f,\kappa) \\
	\le \limsup_{t\to\infty}\limsup_{s\to\infty}e_t^{\mu P_s}(f,\kappa)
	= \limsup_{t\to\infty}e_t(f,\kappa) = e(f,\kappa).
\end{multline*}
This proves lemma \ref{lem:erg1}.  For lemma \ref{lem:erg2}, note that
\begin{equation*}
\begin{split}
	e_{t+s}^\mu(f,\kappa) &= 
	\mathbf{E}^\mu(f(X_{t+s})^2)-
	\mathbf{E}^\mu\!\left(
	\mathbf{E}^\mu\bigg(
	\mathbf{E}^\mu\big(f(X_{t+s})\big|\mathcal{F}^{Y,\kappa}_{[0,t+s]}
	\vee\mathcal{F}^{X}_{[0,s]}\big)
	\bigg|\mathcal{F}^{Y,\kappa}_{[0,t+s]}\bigg)^2
	\right) \\
	&\ge 
	\mathbf{E}^\mu(f(X_{t+s})^2)-
	\mathbf{E}^\mu\!\left(
	\mathbf{E}^\mu
	\big(f(X_{t+s})\big|\mathcal{F}^{Y,\kappa}_{[0,t+s]}
	\vee\mathcal{F}^{X}_{[0,s]}
	\big)^2\right) \\
	&=
	\mathbf{E}^\mu(f(X_{t+s})^2)-
	\mathbf{E}^\mu\!\left(
	\mathbf{E}^\mu\!\left(\left.
	\mathbf{E}^\mu
	\big(f(X_{t+s})\big|\mathcal{F}^{Y,\kappa}_{[s,t+s]}
	\vee\sigma(X_s)
	\big)^2\right|\sigma(X_s)\right)\right) \\
	&=
	\mathbf{E}^{\mu P_s}\left[
	\mathbf{E}^{\delta_{X_0}}(f(X_t)^2)-
	\mathbf{E}^{\delta_{X_0}}\!\left(
	\mathbf{E}^{\delta_{X_0}}\big(f(X_t)\big|\mathcal{F}^{Y,\kappa}_{[0,t]}
	\big)^2\right)
	\right]
	= \mathbf{E}^{\mu P_s}\Big[e_t^{\delta_{X_0}}(f,\kappa)\Big].
\end{split}
\end{equation*}
Thus evidently, we can estimate
$$
	\liminf_{t\to\infty}e_t^\mu(f,\kappa) =
	\liminf_{t\to\infty}\liminf_{s\to\infty}e_{t+s}^\mu(f,\kappa) \ge
	\liminf_{t\to\infty}\mathbf{E}\left[e_t^{\delta_{X_0}}(f,\kappa)\right],
$$
and it remains to establish that the latter limit equals $e(f,\kappa)$.  But
\begin{multline*}
	\left|\mathbf{E}\left[e_t^{\delta_{X_0}}(f,\kappa)\right] -
	e_t(f,\kappa)\right| = 
	\left|
	\mathbf{E}\left[
	\mathbf{E}\big(f(X_t)\big|\mathcal{F}^{Y,\kappa}_{[0,t]}
        \big)^2-
	\mathbf{E}^{\delta_{X_0}}\big(f(X_t)\big|\mathcal{F}^{Y,\kappa}_{[0,t]}
        \big)^2
	\right]\right| \\ \le
	2\|f\|_\infty\,\mathbf{E}\left[
	\left|\mathbf{E}\big(f(X_t)\big|\mathcal{F}^{Y,\kappa}_{[0,t]}
        \big)-
	\mathbf{E}^{\delta_{X_0}}\big(f(X_t)\big|\mathcal{F}^{Y,\kappa}_{[0,t]}
        \big)\right|
	\right]\xrightarrow{t\to\infty}0
\end{multline*}
using $\kappa>0$, ergodicity of the signal, and
\cite[theorem 6.6]{Van08}.
\qquad\endproof

We emphasize, in particular, that in the ergodic case the maximal 
achievable accuracy problem is completely equivalent to the stationary 
maximal achievable accuracy problem for any initial measure $\mu$.
We may therefore concentrate on the stationary case without loss of 
generality, which we will do from now on.  Note, however, that our results 
below do not themselves require ergodicity.

\section{A General Result}
\label{sec:main}

The purpose of this section is to prove the following general result, 
which relates the maximal achievable accuracy problem to the various 
systems theoretic notions introduced above.

\begin{theorem}
\label{thm:main}
The following conditions are equivalent:
\begin{remunerate}
\item The filter achieves maximal accuracy.
\item The filtering model is stably invertible.
\item The filtering model is strongly reconstructible.
\end{remunerate}
Moreover, any of these conditions implies the following:
\begin{remunerate}
\setcounter{muni}{3}
\item The filtering model is invertible.
\item The filtering model is reconstructible.
\end{remunerate}
\end{theorem}

It should be noted that often invertibility and reconstructibility are 
much easier to verify than stable invertibility or strong 
reconstructibility.  However, our general result only shows that the 
former are necessary conditions for the filter to achieve maximal 
accuracy.  In the next section, we will show that when the signal state 
space is a finite set, the filter achieves maximal accuracy if and only if 
the model is both invertible and reconstructible.  This will allow us to 
give simple necessary and sufficient conditions which can be verified 
directly in terms of the model coefficients.  In general, however, it need 
not be the case that invertibility and reconstructibility are sufficient 
for the filter to achieve maximal accuracy, see section \ref{sec:cexample} 
for a counterexample.

The remainder of this section is devoted to the proof of theorem 
\ref{thm:main}.

\subsection{Proof of $1\Leftrightarrow 2$}

The key here is that we can characterize precisely the limit of 
$e(f,\kappa)$ as $\kappa\to 0$.

\begin{lemma}
\label{lem:manywieners}
For every test function $f$ with $\int f^2d\pi<\infty$, we have
$$
	e(f,\kappa)\xrightarrow{\kappa\to 0} e(f,0) =
	\mathbf{E}\left(\left\{
	f(X_0)-\mathbf{E}\big(f(X_0)\big|\mathcal{F}^{h(X)}_{]-\infty,0]}
	\big)
	\right\}^2\right) := e(f).
$$
\end{lemma}

\begin{proof}
Let $\kappa_\ell\searrow 0$ as $\ell\to\infty$.  It suffices to show that
$e(f,\kappa_\ell)\to e(f)$ as $\ell\to\infty$ for every such sequence.  We 
will therefore fix an arbitrary sequence $\kappa_\ell\searrow 0$ in the 
remainder of the proof.

Without loss of generality, we assume (enlarging the probability space if 
necessary) that $(\Omega,\mathcal{F},\mathbb{P})$ carries a countable 
sequence $B_t^\ell$ of independent $n$-dimensional Wiener processes which 
are independent of $(X_t)_{t\in\mathbb{R}}$.  Define
$$
	W_t^r = \sum_{\ell=r}^\infty 
	\sqrt{\kappa_\ell^2-\kappa_{\ell+1}^2}
	\, B_t^\ell,\qquad\quad
	Z_t^r = \int_0^t h(X_s)\,ds + W_t^r.
$$
Note that it is easily established that the sum in the expression for 
$W_t^r$ is a.s.\ convergent uniformly on compact time intervals, and that 
the limit is a Wiener process with covariance $\kappa^2_r\,I$.
The process $(X_t,Y_t^{\kappa_r})_{t\in\mathbb{R}}$ therefore has the same
law as $(X_t,Z_t^r)_{t\in\mathbb{R}}$ under $\mathbf{P}$, and in 
particular (here $\mathcal{F}^{Z,r}_{]-\infty,0]}=\sigma\{Z_t^r:t\le 0\}$)
$$
	e(f,\kappa_r) =
	\mathbf{E}\left(\left\{
	f(X_0)-\mathbf{E}\big(f(X_0)\big|\mathcal{F}^{Z,r}_{]-\infty,0]}
	\big)
	\right\}^2\right).
$$
But by the independence of $B_t^\ell$ and the signal, evidently
$$
	\mathbf{E}\big(f(X_0)\big|\mathcal{F}^{Z,r}_{]-\infty,0]}\big) =
	\mathbf{E}\big(f(X_0)\big|\mathcal{F}^{Z,0}_{]-\infty,0]}
	\vee \mathcal{F}^{B,0}_{]-\infty,0]}\vee\cdots\vee
	\mathcal{F}^{B,r-1}_{]-\infty,0]}\big)\quad
	\mathbf{P}\mbox{-a.s.},
$$
where $\mathcal{F}^{B,\ell}_{]-\infty,0]}=\sigma\{B_t^\ell:t\le 0\}$.
Therefore 
$$
	\mathbf{E}\big(f(X_0)\big|\mathcal{F}^{Z,r}_{]-\infty,0]}\big)
	\xrightarrow{r\to\infty}
	\mathbf{E}\Bigg(f(X_0)\Bigg|\mathcal{F}^{Z,0}_{]-\infty,0]}\vee
	\bigvee_{\ell\ge 0} 
	\mathcal{F}^{B,\ell}_{]-\infty,0]}\Bigg)\quad
	\mbox{in }L^2(\mathbf{P})
$$
by the martingale convergence theorem.  But using again independence
\begin{multline*}
	\mathbf{E}\Bigg(f(X_0)\Bigg|\mathcal{F}^{Z,0}_{]-\infty,0]}\vee
	\bigvee_{\ell\ge 0} 
	\mathcal{F}^{B,\ell}_{]-\infty,0]}\Bigg) = \\
	\mathbf{E}\Bigg(f(X_0)\Bigg|\mathcal{F}^{h(X)}_{]-\infty,0]}\vee
	\bigvee_{\ell\ge 0} 
	\mathcal{F}^{B,\ell}_{]-\infty,0]}\Bigg) =
	\mathbf{E}\big(f(X_0)\big|\mathcal{F}^{h(X)}_{]-\infty,0]}\big)
	\quad\mathbf{P}\mbox{-a.s.},
\end{multline*}
and the claim follows directly.
\qquad
\end{proof}

We can now prove the implications $1\Rightarrow 2$ and $2\Rightarrow 1$ of 
theorem \ref{thm:main}.

{\em Proof of theorem \ref{thm:main}, $1\Rightarrow 2$}.
Suppose the filter achieves maximal accuracy.  Then
$$
	e(f)=0\quad\Longrightarrow\quad
	f(X_0) =
	\mathbf{E}\big(f(X_0)\big|\mathcal{F}^{h(X)}_{]-\infty,0]}\big)
	\quad\mathbf{P}\mbox{-a.s.}
$$
whenever $f$ is bounded and measurable.
As the signal state space $E$ is Polish, it is isomorphic as a measure 
space to a subset of the interval $[0,1]$.  Denote by $\iota:E\to[0,1]$ 
this isomorphism.  Setting $f=\iota$ above, we find that $\iota(X_0)$ 
coincides $\mathbf{P}$-a.s.\ with an 
$\mathcal{F}^{h(X)}_{]-\infty,0]}$-measurable random variable.
Therefore so does $X_0$.
\qquad\endproof

{\em Proof of theorem \ref{thm:main}, $2\Rightarrow 1$}.
Suppose the filtering model is stably invertible.  It follows immediately 
that $e(f)=0$ for every test function $f$ with $\int f^2d\pi<\infty$.
\qquad\endproof

\subsection{Proof of $2\Leftrightarrow 3$}

~\par

{\em Proof of theorem \ref{thm:main}, $2\Rightarrow 3$}.
Let $\mu,\nu\ll\pi$ and $\mu\perp\nu$.  Define the event
$M=\{d\mu/d\pi(X_0)>0\}$.  If the filtering model is stably invertible,
then $I_M$ coincides $\mathbf{P}$-a.s.\ with $I_H$ for some 
$H\in\mathcal{F}^{h(X)}_{]-\infty,0]}$.  But then $\mathbf{P}^\mu(H)=1$ 
and $\mathbf{P}^\nu(H)=0$, so the filtering model is strongly reconstructible.
\qquad\endproof

{\em Proof of theorem \ref{thm:main}, $3\Rightarrow 2$}.
We suppose the model is strongly reconstructible.  Let 
$\{A_1,\ldots,A_m\}$ be a partition of $E$ with $\pi(A_i)>0$ for all $i$.
Define $\pi_i(B)=\pi(B\cap A_i)/\pi(A_i)$.  Then $\pi_i\ll\pi$ for every 
$i$ and $\pi_i\perp\pi_j$ for $i\ne j$.  By strong reconstructibility, we 
may therefore find disjoint $H_1,\ldots,H_m\in\mathcal{F}^{h(X)}_{]-\infty,0]}$
such that $\mathbf{P}^{\pi_i}(H_j)=\delta_{ij}$ for all $i,j$, or, in 
other words, $\mathbf{P}(H_i|X_0\in A_i)=1$.

Now let $f(x)=\sum_if_iI_{A_i}(x)$, and let $H=\sum_if_iI_{H_i}$
($f_i$ are distinct).  Then
\begin{multline*}
	\mathbf{P}(f(X_0)=H) = 
	\sum_{i=1}^m\mathbf{P}(f(X_0)=H|X_0\in A_i)\,\mathbf{P}(X_0\in A_i)
	\\
	= 
	\sum_{i=1}^m\mathbf{P}(H_i|X_0\in A_i)\,\mathbf{P}(X_0\in A_i)
	= 1.
\end{multline*}
Therefore $f(X_0)$ coincides $\mathbf{P}$-a.s.\ with the
$\mathcal{F}^{h(X)}_{]-\infty,0]}$-measurable random variable $H$.

Evidently $f(X_0)$ coincides $\mathbf{P}$-a.s.\ with an
$\mathcal{F}^{h(X)}_{]-\infty,0]}$-measurable random variable whenever $f$ 
is a simple function.  But recall that any measurable function can be 
approximated monotonically by a sequence simple functions, so that the 
claim evidently holds for any measurable function $f$.  It suffices to 
note that as the signal state space $E$ is Polish, it is isomorphic as a 
measure space to a subset of the interval $[0,1]$, so that we may apply 
our conclusion to the isomorphism $\iota$.
\qquad\endproof

\subsection{Proof of $1\Rightarrow 4$}

The proof of this implication follows from the following observation:
it can be read off from the proof of lemmas \ref{lem:erg1} and
\ref{lem:erg2} that 
$$
	\mathbf{E}\left(\left\{
	f(X_t)-\mathbf{E}\big(f(X_t)\big|\mathcal{F}^{h(X)}_{[0,t]}
	\vee\sigma(X_0)\big)
	\right\}^2\right) \le e_{t+s}(f,0)
$$
whenever $\int f^2d\pi<\infty$ and $t,s\ge 0$.
Therefore, if the filter achieves maximal accuracy,
$$
	\mathbf{E}\left(\left\{
	f(X_t)-\mathbf{E}\big(f(X_t)\big|\mathcal{F}^{h(X)}_{[0,t]}
	\vee\sigma(X_0)\big)
	\right\}^2\right) \le \lim_{s\to\infty}e_{t+s}(f,0) =
	e(f)=0.
$$
Thus $f(X_t)$ coincides $\mathbf{P}$-a.s.\ with a
$\mathcal{F}^{h(X)}_{[0,t]}\vee\sigma(X_0)$-measurable random variable.
\qquad\endproof

\subsection{Proof of $3\Rightarrow 5$}

Suppose the model is not reconstructible.  We claim that it cannot be 
strongly reconstructible.  Indeed, if the model is not reconstructible, then 
$$
	\mbox{there exist}\quad\mu,\nu\ll\pi,\quad
	\mu\ne\nu\quad\mbox{such that}\quad
	\mathbf{P}^\mu|_{\mathcal{F}^{h(X)}_{]-\infty,0]}}=
	\mathbf{P}^\nu|_{\mathcal{F}^{h(X)}_{]-\infty,0]}}.
$$
Define $\mu' = (\mu-\nu)^+/(\mu-\nu)^+(E)$ and 
$\nu'=(\mu-\nu)^-/(\mu-\nu)^-(E)$.  Clearly $\mu',\nu'$ are probability 
measures and $\mu'-\nu'\propto \mu-\nu$ (as 
$(\mu-\nu)^+(E)=(\mu-\nu)^-(E)$), so
$$
	\mu',\nu'\ll\pi,\quad
	\mu'\perp\nu',\quad
	\mbox{and}\quad
	\mathbf{P}^{\mu'}|_{\mathcal{F}^{h(X)}_{]-\infty,0]}}=
	\mathbf{P}^{\nu'}|_{\mathcal{F}^{h(X)}_{]-\infty,0]}}.
$$
Hence the model can certainly not be strongly reconstructible.
\qquad\endproof

\section{Finite State Space}
\label{sec:finst}

We have shown above that invertibility and reconstructibility are 
necessary conditions for the filter to achieve maximal accuracy.  In this 
section, we will show that in the case where the signal state space is a 
finite set, these conditions together are also sufficient.  This is 
particularly useful as invertibility and reconstructibility can be 
verified algebraically in terms of the model coefficients, while 
verifying stable invertibility or strong reconstructibility directly is 
difficult.

Let $(X_t)_{t\in\mathbb{R}}$ be a stationary finite state Markov process.
The state space is $E=\{1,\ldots,d\}$, the transition law is determined by 
the $d\times d$ transition intensities matrix $\Lambda=(\lambda_{ij})$, 
and the stationary measure is the $d$-dimensional vector $\pi=(\pi_i)$.  
The observation function is also represented as a $d$-dimensional vector 
$h=(h_i)$ (as no confusion may arise, we will make no distinction between 
functions and measures on $E$ and their representing vectors).  We will 
assume that $\pi_i>0$ for all $i$.

\begin{remark}
The assumption that $\pi_i>0$ for all $i$ is made for convenience only and 
does not entail any loss of generality.  Indeed, if any of the entries of 
the stationary distribution $\pi$ are zero, then we may remove the 
corresponding points from the state space $E$ and apply our results below 
to the resulting stationary Markov process on the reduced state space.
Of course, the algebraic conditions in lemmas \ref{lem:graph} and 
\ref{lem:frecon} below must then be applied to the coefficients of the 
reduced model.
\end{remark}

The main result in this section is the following.

\begin{theorem}
\label{thm:fst}
For the finite state filtering model, the following are equivalent:
\begin{remunerate}
\item The filter achieves maximal accuracy.
\item The filtering model is stably invertible.
\item The filtering model is strongly reconstructible.
\item The filtering model is invertible and reconstructible.
\end{remunerate}
\end{theorem}

Clearly, all that remains to be shown is the implication $4\Rightarrow 1$.  
Before we proceed to the proof, let us show how invertibility and 
reconstructibility can be verified in terms of the model parameters.  To 
this end we give the following two lemmas.

\begin{lemma}
\label{lem:graph}
The finite state filtering model is invertible iff the 
following hold:
\begin{remunerate}
\item For any $i\ne j$ such that $\lambda_{ij}>0$, we have
$h_i\ne h_j$.
\item For any $i\ne j\ne k$ such that $\lambda_{ij}>0$, $\lambda_{ik}>0$,
we have $h_j\ne h_k$.
\end{remunerate}
\end{lemma}

\begin{proof}
Suppose the stated conditions hold.  By the first condition, 
$h(X_t)$ jumps at every jump time $\tau$ of $X_t$.  By the second
condition, if $X_{\tau-}$ is known, then $X_\tau=F[X_{\tau-},h(X_\tau)]$
for a suitable function $F$.  Therefore, if $X_s$ is known for some $s<t$, 
then the entire path $(X_r)_{s\le r\le t}$ can be reconstructed from 
$(h(X_r))_{s\le r\le t}$ by a straightforward algorithm.  This establishes 
invertibility.  Conversely, if one of the stated conditions does not hold, 
it is easily seen that invertibility must fail.
\qquad
\end{proof}

\begin{lemma}
\label{lem:frecon}
The finite state filtering model is reconstructible if and only if
$$
        \mathrm{dim}\left(\mathrm{span}\left\{
        H^{n_0}\tilde\Lambda H^{n_1}\tilde\Lambda
        \cdots\tilde\Lambda H^{n_k}\mathbf{1}:
        k,n_0,\ldots,n_k\ge 0\right\}\right) = d.
$$
Here $\mathbf{1}=(1~~1~~\cdots~~1)^*$ is the column vector of ones; 
$\tilde\Lambda=(\tilde\lambda_{ij})$ 
is the transition intensities matrix whose off-diagonal entries satisfy 
$\tilde\lambda_{ij}=\lambda_{ji}\pi_j/\pi_i$; and $H=\mathrm{diag}(h)$.
\end{lemma}

\begin{proof}
It is readily verified that the time reversed signal $\tilde X_t$ is a 
finite state Markov process with transition intensities matrix 
$\tilde\Lambda$ \cite[theorem 3.7.1]{Nor98}.  As we have assumed without 
loss of generality that every point of the state space is positively 
charged by $\pi$, any probability measure $\mu$ on $E$ is absolutely 
continuous $\mu\ll\pi$.  Therefore reconstructibility in this setting is 
simply observability of the reverse time system, and the condition in the 
lemma follows along the lines of \cite[lemma 9]{Van09}.
\qquad
\end{proof}

The condition in this last lemma can always be computed in a finite number 
of steps; see \cite[remark 11]{Van09} for further comments and a simple 
but explicit algorithm.

\subsection{Proof of $4\Rightarrow 1$}

As we assume invertibility, we have
$$
	f(X_t) = \mathbf{E}\big(f(X_t)\big|\mathcal{F}^{h(X)}_{[0,t]}
        \vee\sigma(X_0)\big)\quad
	\mathbf{P}\mbox{-a.s.}\qquad
	\mbox{for all functions }f.
$$
Therefore, we can evidently write
\begin{equation*}
\begin{split}
	e_t(f,0) &= 
	\mathbf{E}\left(\left\{
	\mathbf{E}\big(f(X_t)\big|\mathcal{F}^{h(X)}_{[0,t]}
        \vee\sigma(X_0)\big)
	-\mathbf{E}\big(f(X_t)\big|\mathcal{F}^{h(X)}_{[0,t]}\big)
	\right\}^2\right) \\
	&=
	\mathbf{E}\left(\left\{
	\mathbf{E}\big(f(\tilde X_0)\big|\mathcal{F}^{h(\tilde X)}_{[0,t]}
        \vee\sigma(\tilde X_{t})\big)
	-\mathbf{E}\big(f(\tilde X_0)\big|\mathcal{F}^{h(\tilde X)}_{[0,t]}\big)
	\right\}^2\right).
\end{split}
\end{equation*}
We would like to show that $e(f)=0$, i.e., that $e_t(f,0)\to 0$ as 
$t\to\infty$ for any $f$.  As 
$$
	e_t(\alpha f+\beta g,0)\le
	2\alpha^2\,e_t(f,0)+2\beta^2\,e_t(g,0),
$$
it clearly suffices to restrict to positive $f>0$ such that $\int 
fd\pi=1$.  Fix such a function, and define the probability measure 
$d\nu=fd\pi$.  Then
\begin{equation*}
\begin{split}
	e_t(f,0) &\le
	2\,\|f\|_\infty\,
	\mathbf{E}\left(\left|
	\mathbf{E}\big(f(\tilde X_0)\big|\mathcal{F}^{h(\tilde X)}_{[0,t]}
        \vee\sigma(\tilde X_{t})\big)
	-\mathbf{E}\big(f(\tilde X_0)\big|\mathcal{F}^{h(\tilde X)}_{[0,t]}\big)
	\right|\right) \\
	&=
	2\,\|f\|_\infty\,
	\mathbf{E}^\nu\left(
	\mathbf{E}\left(\left.\left|\frac{
	\mathbf{E}\big(f(\tilde X_0)\big|\mathcal{F}^{h(\tilde X)}_{[0,t]}
        \vee\sigma(\tilde X_{t})\big)}{
	\mathbf{E}\big(f(\tilde X_0)\big|\mathcal{F}^{h(\tilde X)}_{[0,t]}\big)
	}-1\right|~\right|\mathcal{F}^{h(\tilde X)}_{[0,t]}\right)\right) \\
	&=
	2\,\|f\|_\infty\, 
	\mathbf{E}^\nu\big(
	\big\|
	\mathbf{P}^\nu(\tilde X_t\in\,\cdot\,|\mathcal{F}^{h(\tilde X)}_{[0,t]})
	-\mathbf{P}(\tilde X_t\in\,\cdot\,|\mathcal{F}^{h(\tilde X)}_{[0,t]})
	\big\|_{\rm TV}\big) \\
	&=
	2\,\|f\|_\infty\sum_{i=1}^d
	\mathbf{E}^\nu\big(\big|
	\mathbf{P}^\nu(\tilde X_t=i|\mathcal{F}^{h(\tilde X)}_{[0,t]})
	-\mathbf{P}(\tilde X_t=i|\mathcal{F}^{h(\tilde X)}_{[0,t]})
	\big|\big),
\end{split}
\end{equation*}
where we used the Bayes formula as in \cite[lemma 5.6 and 
corollary 5.7]{Van08}
to show that
\begin{multline*}
	\big\|
	\mathbf{P}^\nu(\tilde X_t\in\,\cdot\,|\mathcal{F}^{h(\tilde X)}_{[0,t]})
	-\mathbf{P}(\tilde X_t\in\,\cdot\,|\mathcal{F}^{h(\tilde X)}_{[0,t]})
	\big\|_{\rm TV} = \\
	\mathbf{E}\left(\left.\left|\frac{
	\mathbf{E}\big(f(\tilde X_0)\big|\mathcal{F}^{h(\tilde X)}_{[0,t]}
        \vee\sigma(\tilde X_{t})\big)}{
	\mathbf{E}\big(f(\tilde X_0)\big|\mathcal{F}^{h(\tilde X)}_{[0,t]}\big)
	}-1\right|~\right|\mathcal{F}^{h(\tilde X)}_{[0,t]}\right).
\end{multline*}
But by reconstructibility and \cite[corollary 1]{Van09}, we find that
$$
	e_t(f,0) \le
	2\,\|f\|_\infty\sum_{i=1}^d
	\mathbf{E}^\nu\big(\big|
	\mathbf{P}^\nu(\tilde X_t=i|\mathcal{F}^{h(\tilde X)}_{[0,t]})
	-\mathbf{P}(\tilde X_t=i|\mathcal{F}^{h(\tilde X)}_{[0,t]})
	\big|\big) \xrightarrow{t\to\infty}0.
$$
Therefore $e_t(f,0)\to 0$ as $t\to\infty$ for any $f$, and the proof is 
complete.
\qquad\endproof

\begin{remark}
\label{rem:gap}
It is interesting to note that all the steps in this proof have 
counterparts in the general setting of section \ref{sec:prelim}.
In particular, it is not difficult to establish that in general, 
to achieve maximal accuracy it is sufficient that the model is invertible 
and that the time-reversed noiseless filter is stable in the sense that
$$
	\mathbf{E}^\nu\big(\big\|
	\mathbf{P}^\nu(\tilde X_t\in\,\cdot\,|\mathcal{F}^{h(\tilde X)}_{[0,t]})
	-\mathbf{P}(\tilde X_t\in\,\cdot\,|\mathcal{F}^{h(\tilde X)}_{[0,t]})
	\big\|_{\rm TV}\big)
	\xrightarrow{t\to\infty}0\quad\forall~
	\nu\ll\pi,~\|\tfrac{d\nu}{d\pi}\|_\infty<\infty.
$$
On the other hand, the results in \cite{Van09} are easily adapted to show 
that if the model is reconstructible, then the time-reversed noiseless 
filter is stable in the sense that 
$$
	\mathbf{E}^\nu\big(\big|
	\mathbf{E}^\nu(g(\tilde X_t)|\mathcal{F}^{h(\tilde X)}_{[0,t]})
	-\mathbf{E}(g(\tilde X_t)|\mathcal{F}^{h(\tilde X)}_{[0,t]})
	\big|\big)
	\xrightarrow{t\to\infty}0
	\quad\forall~\nu\ll\pi,~\|\tfrac{d\nu}{d\pi}\|_\infty<\infty,~
	g\in L^2(\pi).
$$
Therefore, invertibility and reconstructibility imply maximal accuracy if 
one can close the gap between total variation stability and individual
stability of the time-reversed noiseless filter.  This is automatic in a 
finite state space (trivial) and in a countable state space (as the 
sequence space $\ell_1$ has the Schur property \cite[theorem 4.32]{AB06}).  
When the signal state space is continuous, however, invertibility and 
reconstructibility is typically not sufficient to guarantee that the 
filter achieves maximal accuracy; a counterexample is given in the next 
section.
\end{remark}

\section{Linear Gaussian Models}
\label{sec:ks}

In this section, we consider a linear Gaussian model of the following 
form ($t\in\mathbb{R}$):
\begin{equation*}
\begin{split}
	X_t &= X_0 + \int_0^t AX_u\,du + D\,W_t, \\
	Y_t^\kappa &= \int_0^t HX_u\,du + \kappa\,B_t.
\end{split}
\end{equation*}
Here $(B_t)_{t\in\mathbb{R}}$ and $(W_t)_{t\in\mathbb{R}}$ are independent 
two-sided Wiener processes of dimensions $n$ and $m$, respectively,
and the signal state space is $E=\mathbb{R}^p$.  We refer to \cite{KS72bk}
for the definitions and basic properties of the various elements of linear 
systems theory (stable matrix, detectability, stabilizability, etc.)\ 
which are used below.

We make the following assumptions:
\begin{remunerate}
\item $A\in\mathbb{R}^{p\times p}$ is an asymptotically stable matrix
and $(X_t)_{t\in\mathbb{R}}$ is stationary;
\item $D\in\mathbb{R}^{p\times m}$ and $H\in\mathbb{R}^{n\times p}$
are matrices of full rank and $m,n\le p$.
\end{remunerate}
The stability of $A$ ensures that the signal is ergodic, and in particular 
that the stationary solution of the signal equation exists and is unique.  

\begin{remark}
The rank assumption on $D$ and $H$ and the assumption on the dimensions 
$m,n,p$ is made for convenience only and does not entail any loss of 
generality.  Indeed, when the matrices are not of full rank we can 
trivially obtain an equivalent model of full rank by reducing the 
dimensions of $W$, $B$ and/or $Y^\kappa$.  Similarly, if $m>p$ we can
obtain an equivalent model with $m=p$.  Of course, the algebraic 
condition in theorem \ref{thm:KS} must then be applied to the coefficients 
of the reduced model.  If $n>p$ the filter is trivially seen to achieve 
maximal accuracy (as then $H$, being of full rank, has a left inverse, so 
the noiseless observations are fully informative).
\end{remark}

The maximal achievable accuracy problem in the linear Gaussian setting was 
considered in a classic paper of Kwakernaak and Sivan \cite{KS72}, where 
an almost\footnote{
	The result of Kwakernaak and Sivan, stated as corollary 
	\ref{cor:ks} below, gives necessary and sufficient conditions
	for the case $m\ge n$ but only a sufficient condition for the case
	$m<n$.
}
necessary and sufficient condition was obtained.
Their approach is surprisingly complicated, however, and relies on rather 
explicit computations of the behavior of Riccati equations in the limit of 
vanishing noise.  In this section we give a direct proof of their theorem 
by verifying the stable invertibility property.

\begin{theorem}
\label{thm:KS}
In the linear Gaussian setting of this section, the filter achieves 
maximal accuracy if and only if the matrix $H(\lambda I-A)^{-1}D$ has 
linearly independent columns for all $\lambda\in\mathbb{C}$ with 
$\mathrm{Re}\,\lambda>0$.
\end{theorem}

The result of Kwakernaak and Sivan follows easily.

\begin{corollary}[Kwakernaak-Sivan]
\label{cor:ks}
The following hold.
\begin{remunerate}
\item If $m>n$, then the filter does not achieve maximal accuracy.
\item If $m=n$, then the filter achieves maximal accuracy if and only
if $\det[H(\lambda I-A)^{-1}D]$ is nonzero for any $\lambda\in\mathbb{C}$ 
with $\mathrm{Re}\,\lambda>0$.
\item If $m<n$ and there exists $M\in\mathbb{R}^{m\times n}$ such that
$\det[MH(\lambda I-A)^{-1}D]$ is nonzero for any $\lambda\in\mathbb{C}$
with $\mathrm{Re}\,\lambda>0$, then the filter achieves maximal accuracy.
\end{remunerate}
\end{corollary}

\begin{proof}
We consider each case separately.
\begin{remunerate}
\item $H(\lambda I-A)^{-1}D$ has $m$ columns each of which is an
$n$-dimensional vector.  Therefore if $m>n$, the columns cannot be 
linearly independent for any $\lambda$.
\item When $m=n$ the matrix $H(\lambda I-A)^{-1}D$ is square, so that it 
has linearly independent columns if and only if 
$\det[H(\lambda I-A)^{-1}D]\ne 0$.
\item If $\det[MH(\lambda I-A)^{-1}D]$ is nonzero, the square matrix
$MH(\lambda I-A)^{-1}D$ has linearly independent columns.  Then 
certainly the columns of $H(\lambda I-A)^{-1}D$ are linearly independent.
\end{remunerate}
In view of these facts, the corollary follows by applying theorem 
\ref{thm:KS}.
\qquad
\end{proof}

\begin{remark} As is pointed out by Godbole \cite{God72}, the condition of 
theorem \ref{thm:main} corresponds to the requirement that the model is 
invertible and that the inverse has no unstable modes (in the sense of 
linear systems).  Indeed, note that $H(\lambda I-A)^{-1}D$ is the transfer 
function associated to our filtering model, so that invertibility holds in 
this setting if and only if the matrix $H(\lambda I-A)^{-1}D$ has a left 
inverse for all but a finite number of $\lambda\in\mathbb{C}$ (see, e.g., 
\cite[theorem 5]{Moy77}).  The inverse is again a linear system whose 
transfer function is the left inverse of $H(\lambda I-A)^{-1}D$, so that 
the lack of right halfplane zeros of $H(\lambda I-A)^{-1}D$ ensures that 
the inverse system does not have any unstable poles.  If there are 
additionally no zeros on the imaginary axis, then the inverse system is 
even asymptotically stable and the heuristic outlined in remark
\ref{rem:stabinv} can be rigorously implemented.

However, it is not immediately obvious from such arguments that stable 
invertibility follows even when there are zeros on the imaginary axis, or 
that the model cannot be stably invertible when there are right halfplane 
zeros.  In the proof of theorem \ref{thm:KS}, the former problem is 
circumvented by using the idea in \cite{KS72} of using an approximate, 
rather than exact, inverse system.  The latter problem is easily resolved
directly in our setting, and the proof of this part of the theorem is 
substantially simpler than the corresponding arguments of Kwakernaak and 
Sivan.
\end{remark}

\subsection{An example}
\label{sec:cexample}

Before we proceed to the proof of theorem \ref{thm:KS}, let us demonstrate 
by means of an example that, unlike in the finite state setting, 
invertibility and reconstructibility are not always sufficient to ensure 
that the filter achieves maximal accuracy.  This implies the existence of 
a gap between total variation and individual stability of the time 
reversed filter, discussed in remark \ref{rem:gap}.

For our example, let $m=n=1$ and $p=2$, and we set
$$
        A = \left[\begin{array}{rr} -1 & 0 \\ 0 & -4 \end{array}\right],
        \qquad
        D = \left[\begin{matrix}1 \\ 1\end{matrix}\right],\qquad
        H = \left[1\quad{-2}\right].
$$
This model satisfies all the assumptions of this section.
We now compute
$$
        H(\lambda I-A)^{-1}D =
        \left[1\quad{-2}\right] \left[\begin{array}{cc} (\lambda+1)^{-1} 
        & 0 \\ 0 & 
        (\lambda+4)^{-1}\end{array}\right]\left[\begin{matrix}1 \\ 
        1\end{matrix}\right] =
        (\lambda+1)^{-1}-2(\lambda+4)^{-1}.
$$
Therefore $H(\lambda I-A)^{-1}D = 0$ for $\lambda=2$, so by theorem 
\ref{thm:KS} the filter does not achieve maximal accuracy.  However, we 
claim that the model is both invertible and reconstructible.

To prove invertibility, it suffices to note that $H(\lambda I-A)^{-1}D$
is nonzero (hence left invertible) for all but a finite number of 
$\lambda\in\mathbb{C}$.  To prove reconstructibility, we use the fact that 
the reverse time signal $\tilde X_t$ satisfies an equation of the form
\cite{And82}
$$
        \tilde X_t = \tilde X_0 + \int_0^t
	\Sigma A^*\Sigma^{-1}\tilde X_s\,ds + D\,\tilde W_t,
$$
where $\tilde W_t$ is a suitably defined Wiener process and $\Sigma$ 
denotes the stationary covariance matrix of the signal.  The matrix
$\Sigma$ can be computed as the unique solution of the Lyapunov equation:
$$
	A\Sigma + \Sigma A^* + DD^* = 0\qquad
	\Longrightarrow\qquad
        \Sigma = \left[\begin{array}{rr} 1/2 & 1/5 \\ 1/5 & 1/8 
	\end{array}\right].
$$
Note that $\Sigma$ is a strictly positive matrix.  Therefore, the model is 
evidently reconstructible if $H\Sigma$ and $H\Sigma A^*$ are linearly 
independent (so that the time reversed model is observable).  But this is 
easily established to be the case by explicit computation.

\subsection{Proof of Theorem \ref{thm:KS}}

We will show that the condition of the theorem is necessary and sufficient 
for stable invertibility.  The necessity part of the proof is closely 
related to a problem of Karhunen \cite{Kar49}, while sufficiency is proved 
along the lines of \cite{KS72}.

In the following, let us write $Z_t=HX_t$ ($t\in\mathbb{R}$) for 
notational simplicity.  We introduce the Hilbert spaces of 
random variables $\mathcal{L}_X,\mathcal{L}_Z,\mathcal{L}_W\subset 
L^2(\mathbf{P})$ as follows:
\begin{equation*}
\begin{split}
	\mathcal{L}_X &= L^2(\mathbf{P})\mbox{-}\mathrm{cl}\,\{
	v_1^*X_{t_1}+\cdots+v_k^*X_{t_k}:k\in\mathbb{N},~
		t_1,\ldots,t_k\le 0,~v_1,\ldots,v_k\in\mathbb{R}^p\},\\
	\mathcal{L}_Z &= L^2(\mathbf{P})\mbox{-}\mathrm{cl}\,\{
	v_1^*Z_{t_1}+\cdots+v_k^*Z_{t_k}:k\in\mathbb{N},~
		t_1,\ldots,t_k\le 0,~v_1,\ldots,v_k\in\mathbb{R}^n\},\\
	\mathcal{L}_W &= L^2(\mathbf{P})\mbox{-}\mathrm{cl}\,\{
	v_1^*W_{t_1}+\cdots+v_k^*W_{t_k}:k\in\mathbb{N},~
		t_1,\ldots,t_k\le 0,~v_1,\ldots,v_k\in\mathbb{R}^m\}.
\end{split}
\end{equation*}
For an $\mathbb{R}^k$-valued random variable $K$ we will write, e.g.,
$K\in\mathcal{L}_X$ when $v^*K\in\mathcal{L}_X$ for every 
$v\in\mathbb{R}^k$.

As the joint process $(X_t,Z_t,W_t)_{t\in\mathbb{R}}$ is Gaussian, the 
stable invertibility problem is essentially linear and can be reduced to 
the investigation of the spaces $\mathcal{L}_X,\mathcal{L}_Z,\mathcal{L}_W$.

\begin{lemma}
The model is stably invertible if and only if 
$\mathcal{L}_Z=\mathcal{L}_W$.
\end{lemma}

\begin{proof}
By definition, the model is stably invertible iff $X_0$ coincides 
$\mathbf{P}$-a.s.\ with a $\sigma\{Z_s:s\le 0\}$-measurable random 
variable, i.e., iff $\mathbf{E}(X_0|Z_{]-\infty,0]})=X_0$.  However, as 
$(X_0,Z_s:s\le 0)$ is Gaussian, it is well known that 
$\mathbf{E}(X_0|Z_{]-\infty,0]})\in\mathcal{L}_Z$ (see, e.g.,
\cite[lemma 6.2.2]{Oks98}).  The model is therefore 
stably invertible iff $X_0\in\mathcal{L}_Z$.

To proceed, note that
$$
	X_t = \int_{-\infty}^t e^{A(t-s)}D\,dW_s,\qquad
	\quad
	Z_t = \int_{-\infty}^t He^{A(t-s)}D\,dW_s\qquad
	\quad (t\in\mathbb{R})
$$
as $A$ is asymptotically stable.  Therefore clearly $X_0\in\mathcal{L}_W$ 
and $\mathcal{L}_Z\subset\mathcal{L}_W$.  If $\mathcal{L}_W=\mathcal{L}_Z$, 
it then follows immediately that $X_0\in\mathcal{L}_Z$.  
It therefore remains to show that $X_0\in\mathcal{L}_Z$ implies
$\mathcal{L}_W\subset\mathcal{L}_Z$.  To this end, assume 
$X_0\in\mathcal{L}_Z$.  By stationarity, $X_t\in\mathcal{L}_Z$ also 
for $t\le 0$.  Therefore evidently $\mathcal{L}_X\subset\mathcal{L}_Z$.  
But
$$
	DW_s = X_s - X_0 + \int_s^0 AX_u\,du,\qquad
	s\le 0.
$$
As $D$ has full rank, we find that $\mathcal{L}_W\subset\mathcal{L}_X$.  
Therefore $\mathcal{L}_W\subset\mathcal{L}_Z$ as required.
\qquad
\end{proof}

We can now complete the proof of theorem \ref{thm:KS}.

{\em Proof of theorem \ref{thm:KS}}.
Suppose $H(\lambda I-A)^{-1}D$ does not have linearly independent columns 
for some $\lambda\in\mathbb{C}$ with $\mathrm{Re}\,\lambda>0$.  Then there 
exists $0\ne w\in\mathbb{C}^m$ such that $H(\lambda I-A)^{-1}Dw=0$, and we 
define
$$
	U = U_1+iU_2 := \int_{-\infty}^0 (e^{\lambda s}w)^*dW_s,\qquad
	U_1,U_2\in\mathcal{L}_W.
$$
We can now compute
$$
	\mathbf{E}(U^*v^*Z_u) 
	= 
	\int_{-\infty}^u e^{\lambda s}\,v^*He^{A(u-s)}Dw\,ds
	=
	e^{\lambda u}\,v^*H(\lambda I-A)^{-1}Dw = 0
$$
for all $u\le 0$, $v\in\mathbb{R}^n$.
In particular, as $v^*Z_u$ is real-valued, 
$\langle{U_1,v^*Z_u}\rangle_{L^2(\mathbf{P})} =
\langle{U_2,v^*Z_u}\rangle_{L^2(\mathbf{P})} = 0$
for $u\le 0$, $v\in\mathbb{R}^n$, so that $U_1,U_2\perp\mathcal{L}_Z$.  
But as $U$ is nonzero, evidently $\mathcal{L}_Z\ne\mathcal{L}_W$ and the 
model is not stably invertible.

Conversely, suppose that the matrix $H(\lambda I-A)^{-1}D$ has
linearly independent columns for all $\lambda\in\mathbb{C}$ with
$\mathrm{Re}\,\lambda>0$.  We will prove that for any $\varepsilon>0$,
there is a random variable $X_0^\varepsilon$ of the form
$$
	X_0^\varepsilon = \int_{-\infty}^0 m_\varepsilon(s)\,Z_s\,ds
	\in \mathcal{L}_Z
$$
such that $\|X_0^\varepsilon-X_0\|_{L^2(\mathbf{P})}<\varepsilon$.
Then $X_0\in\mathcal{L}_Z$ (as $\mathcal{L}_Z$ is closed),
so the model is stably invertible.

To prove the claim, fix $\varepsilon>0$.  Then
$$
	X_0^\varepsilon = 
	\int_{-\infty}^0 
	m_\varepsilon(s)
	\int_{-\infty}^0 I_{s\ge u} He^{A(s-u)}D\,dW_u\,ds
	=
	\int_{-\infty}^0 
	\int_{u}^0 m_\varepsilon(s)\,He^{A(s-u)}D\,ds\,dW_u,
$$
provided that $m_\varepsilon$ is bounded and the function
$$
	T_\varepsilon(u) := \int_{u}^0 m_\varepsilon(s)\,He^{A(s-u)}D\,ds
$$
is square integrable (this is justified by truncating the lower bounds 
on the integrals and applying Fubini's theorem for stochastic integrals
\cite[theorem IV.64]{Pro04}).  Note that
$$
	\|X_0^\varepsilon-X_0\|_{L^2(\mathbf{P})}^2 = 
	\int_{-\infty}^0 \|T_\varepsilon(u)-e^{-Au}D\|_F^2\,du,
$$
where $\|C\|_F^2=\mathrm{Tr}[CC^*]$ is the Frobenius norm.  Define for 
$\mathrm{Re}\,\lambda>0$ the Laplace transforms
$$
	\hat m_\varepsilon(\lambda) = \int_{-\infty}^0 e^{\lambda s}\,
	m_\varepsilon(s)\,ds,
	\qquad
	\hat T_\varepsilon(\lambda) = \int_{-\infty}^0 e^{\lambda s}\,
	T_\varepsilon(s)\,ds = 
	\hat m_\varepsilon(\lambda)\,H(\lambda I-A)^{-1}D.
$$
By Plancherel's theorem, we can write (see, e.g., \cite[pp.\ 162--163]{Yos80})
\begin{multline*}
	\|X_0^\varepsilon-X_0\|_{L^2(\mathbf{P})}^2 = \\
	\frac{1}{2\pi}\lim_{x\searrow 0}\int_{-\infty}^\infty
	\|\hat m_\varepsilon(x+iy)\,H(\{x+iy\}I-A)^{-1}D
	-(\{x+iy\}I-A)^{-1}D\|_F^2\,dy.
\end{multline*}
It remains to choose $m_\varepsilon$ with the required properties
such that this expression is smaller than $\varepsilon^2$.  

By our assumption, the left inverse $V(\lambda)$ of the matrix 
$H(\lambda I-A)^{-1}D$ is defined on the right halfplane, i.e., 
$V(\lambda)H(\lambda I-A)^{-1}D = I$ for $\mathrm{Re}\,\lambda>0$.
The above expression for $\|X_0^\varepsilon-X_0\|_{L^2(\mathbf{P})}$
is therefore identically zero if we choose $\hat m_\varepsilon(\lambda)$ 
as $(\lambda I-A)^{-1}DV(\lambda)$.  The problem is that the latter may 
not be the Laplace transform of a function $m_\varepsilon$ with the 
required properties.  We therefore regularize as follows:
$$
	\hat m_\varepsilon(\lambda) = 
	\frac{\gamma^\ell}{(\lambda+\gamma)^\ell}\,
	(\lambda I-A)^{-1}D\,V(\lambda)
$$
for some $\gamma>0$, $\ell\in\mathbb{N}$ to be chosen presently.  As 
$\lambda\mapsto H(\lambda I-A)^{-1}D$ is a rational function, 
$\hat m_\varepsilon$ is rational also.  We choose $\ell$ sufficiently 
large that the degree of the denominator is larger than the degree of the 
numerator.  Then $\hat m_\varepsilon$ is strictly proper with poles in the 
closed left halfplane $\mathrm{Re}\,\lambda\le 0$ only, and is therefore 
the Laplace transform of some bounded function $m_\varepsilon$.  Moreover, 
as
\begin{multline*}
	\sup_{x>0}\int_{-\infty}^\infty
	\|\hat T_\varepsilon(x+iy)\|_F^2\,dy = \\
	\sup_{x>0}\int_{-\infty}^\infty
	\left|\frac{\gamma^\ell}{(x+iy+\gamma)^\ell}\right|^2
	\|(\{x+iy\}I-A)^{-1}D\|_F^2\,dy \le
	2\pi\,\|X_0\|_{L^2(\mathbf{P})}^2,
\end{multline*}
the function $T_\varepsilon$ is square integrable by the Paley-Wiener 
theorem.  Finally, as
$$
	\|X_0^\varepsilon-X_0\|_{L^2(\mathbf{P})}^2 =
	\frac{1}{2\pi}\int_{-\infty}^\infty
	\left|\frac{\gamma^\ell}{(iy+\gamma)^\ell}-1\right|^2
	\|(iyI-A)^{-1}D\|_F^2\,dy \xrightarrow{\gamma\to\infty}0
$$
by dominated convergence, we may choose $\gamma$ such that
$\|X_0^\varepsilon-X_0\|_{L^2(\mathbf{P})}<\varepsilon$.
\qquad\endproof

\subsection{The unstable case}

To be fair, it should be noted that we have not entirely reproduced the 
result of Kwakernaak and Sivan as we have assumed that the matrix $A$ is 
asymptotically stable.  The result of Kwakernaak and Sivan states that the 
conclusion of corollary \ref{cor:ks} above holds already under the weaker 
assumption that the filtering model is detectable and stabilizable.  
Unfortunately, our approach relies crucially on the stationarity (or 
ergodicity) of the signal process, so that one could never expect to 
obtain general results in the setting where the signal may be transient.  
On the other hand, in the linear Gaussian case, the special structure of 
the model allows us to reduce the detectable/stabilizable case to the 
stationary case.  We therefore recover the result of Kwakernaak and Sivan 
in its entirety.

We develop the relevant argument presently.  Let us emphasize, however, 
that the following argument is very specific to the linear Gaussian setting.

We consider again a linear Gaussian model of the form
\begin{equation*}
\begin{split}
	X_t &= X_0 + \int_0^t AX_u\,du + D\,W_t, \\
	Y_t^\kappa &= \int_0^t HX_u\,du + \kappa\,B_t,
\end{split}
\end{equation*}
where $(B_t)_{t\in\mathbb{R}}$ and $(W_t)_{t\in\mathbb{R}}$ are 
independent two-sided Wiener processes of dimensions $n$ and $m$, 
respectively, and the signal state space is $E=\mathbb{R}^p$. 
As above, we will assume $D\in\mathbb{R}^{p\times m}$ and 
$H\in\mathbb{R}^{n\times p}$ are matrices of full rank and $m,n\le p$.
We do not assume, however, that $A\in\mathbb{R}^{p\times p}$ is stable;
instead, we assume only that $(A,D)$ is stabilizable 
and $(A,H)$ is detectable.  The law of $X_0$ may be chosen arbitrarily.

As $(A,H)$ is detectable, it is well known that there exists a matrix 
$K\in\mathbb{R}^{p\times n}$ such that $\bar A:=A-KH$ is an asymptotically 
stable matrix.  Fix such a matrix $K$ (it may not be unique, but this will 
not affect our final result).  By It\^o's rule,
\begin{multline*}
	e^{-\bar At}X_t = 
	X_0 + \int_0^t e^{-\bar As}KHX_s\,ds + \int_0^t e^{-\bar As}D\,dW_s = 
	\\
	X_0 + \int_0^t e^{-\bar As}K\,dY_s^\kappa 
	- \kappa \int_0^t e^{-\bar As}K\,dB_s
	+ \int_0^t e^{-\bar As}D\,dW_s.
\end{multline*}
Now define
$$
	\bar X_t^\kappa := X_t - \int_0^t e^{\bar A(t-s)}K\,dY_s^\kappa,
	\qquad\quad
	\bar Y_t^\kappa := \int_0^t H\bar X_u^\kappa\,du + \kappa\,B_t.
$$
Then evidently $\bar X_t^\kappa$ satisfies the stochastic differential 
equation
$$
	\bar X_t^\kappa = X_0 + \int_0^t \bar A\bar X_s^\kappa\,ds
	+ D\,W_t - \kappa\,K\,B_t.
$$
Moreover, we can compute
\begin{multline*}
	e^{-At}\bar X_t^\kappa = 
	X_0 + \int_0^t e^{-As}D\,dW_s 
	- \int_0^t e^{-As}KH\bar X_s^\kappa\,ds
	- \kappa\int_0^t e^{-As}K\,dB_s = \\
	X_0 + \int_0^t e^{-As}D\,dW_s - \int_0^t e^{-As}K\,d\bar Y_s^\kappa.
\end{multline*}
Thus evidently
$$
	X_t = \bar X_t^\kappa + \int_0^t e^{A(t-s)}K\,d\bar Y_s^\kappa.
$$
The following lemma is therefore immediate.

\begin{lemma}
$\sigma\{Y_t^\kappa:t\in[0,T]\} = \sigma\{\bar Y_t^\kappa:t\in [0,T]\}$ 
$\mathbf{P}$-a.s.\ $\forall$ $T\le\infty$, $\kappa\ge 0$.
\end{lemma}

\begin{proof}
This follows directly from
$$
	\bar Y_t^\kappa = Y_t^\kappa
	- \int_0^t\int_0^s He^{\bar A(s-u)}K\,dY_u^\kappa\,ds,\qquad
	Y_t^\kappa = \bar Y_t^\kappa 
	+ \int_0^t\int_0^s He^{A(s-u)}K\,d\bar Y_u^\kappa\,ds.
$$
The proof is complete.
\qquad
\end{proof}

We now see immediately that for every $t,\kappa\ge 0$
$$
	\mathbf{E}\left(
	\left\|X_t - \mathbf{E}(X_t|\mathcal{F}_{[0,t]}^{Y,\kappa})\right\|^2
	\right) = \mathbf{E}\left(\left\|\bar X_t^\kappa - 
	\mathbf{E}(\bar X_t^\kappa|\mathcal{F}_{[0,t]}^{\bar Y,\kappa})\right\|^2
	\right),
$$
where $\mathcal{F}_{[0,t]}^{\bar Y,\kappa}$ is defined in the obvious 
fashion.  But $\bar X_t^\kappa$ is an ergodic Markov process (as $\bar A$ 
is stable), which brings us back---in principle---to the setting 
employed throughout this paper.  However, note that the driving noise of 
$\bar X_t^\kappa$ is correlated with the observation noise, so that we 
can not immediately apply our previous results.

\begin{lemma}
Define for $t\in\mathbb{R}$
$$
	X_t^0 := \int_{-\infty}^t e^{\bar A(t-s)}D\,dW_s.
$$
Then we have
$$
	\lim_{\kappa\to 0}\lim_{t\to\infty}\mathbf{E}\left(
	\left\|X_t - \mathbf{E}(X_t|\mathcal{F}_{[0,t]}^{Y,\kappa})\right\|^2
	\right) =
	\mathbf{E}\left(
	\left\|X_0^0 - \mathbf{E}(X_0^0|
	\sigma\{HX_s^0:s\le 0\})\right\|^2
	\right).
$$
\end{lemma}

\begin{proof}
Let us define, for $\kappa\ge 0$ and $t\ge 0$, the processes
$$
	x_t^\kappa = x_0^\kappa + \int_0^t Ax_u^\kappa\,du + D\,W_t,
	\qquad
	x_0^\kappa = \int_{-\infty}^0 e^{-\bar As}D\,dW_s - 
	\kappa\int_{-\infty}^0 e^{-\bar As}K\,dB_s,
$$
$$
	y_t^\kappa = \int_0^t Hx_u^\kappa\,du + \kappa\,B_t,\qquad
	\bar x_t^\kappa = x_t^\kappa 
	- \int_0^t e^{\bar A(t-s)}K\,dy_s^\kappa,\qquad
	\bar y_t^\kappa = \int_0^t H\bar x_u^\kappa\,du + \kappa B_t.
$$
Then $(x_t^\kappa,y_t^\kappa)$ is a Markov process with the same 
transition law as $(X_t,Y_t^\kappa)$, except that we have chosen a 
specific initial law for $x_0^\kappa$ in a manner that depends on 
$\kappa$.  However, as the model is stabiliziable and detectable, it is 
well known that the stationary filtering error exists and is independent 
of the initial law.  Therefore
$$
	\lim_{t\to\infty}\mathbf{E}\left(
	\left\|X_t - \mathbf{E}(X_t|\mathcal{F}_{[0,t]}^{Y,\kappa})\right\|^2
	\right) =
	\lim_{t\to\infty}\mathbf{E}\left(
	\left\|x_t^\kappa - \mathbf{E}(x_t^\kappa|\mathcal{F}_{[0,t]}^{y,\kappa})\right\|^2
	\right).
$$
From the above discussion, it is now easily seen that in fact also
$$
	\lim_{t\to\infty}\mathbf{E}\left(
	\left\|X_t - \mathbf{E}(X_t|\mathcal{F}_{[0,t]}^{Y,\kappa})\right\|^2
	\right) =
	\lim_{t\to\infty}\mathbf{E}\left(
	\left\|\bar x_t^\kappa - \mathbf{E}(\bar x_t^\kappa|\mathcal{F}_{[0,t]}^{\bar y,\kappa})\right\|^2
	\right).
$$
Here we have defined $\mathcal{F}_{[0,t]}^{y,\kappa}$ and 
$\mathcal{F}_{[0,t]}^{\bar y,\kappa}$ in the obvious fashion.

Now note that $\bar x_t^\kappa$ is a stationary Markov process with the 
explicit representation
$$
	\bar x_t^\kappa =
	\int_{-\infty}^t e^{\bar A(t-s)}D\,dW_s - 
	\kappa\int_{-\infty}^t e^{\bar A(t-s)}K\,dB_s.
$$
Therefore $(x_t^\kappa,y_t^\kappa)$ immediately extend to all 
$t\in\mathbb{R}$, and by stationarity
$$
	\lim_{t\to\infty}\mathbf{E}\left(
	\left\|X_t - \mathbf{E}(X_t|\mathcal{F}_{[0,t]}^{Y,\kappa})\right\|^2
	\right) =
	\mathbf{E}\left(
	\left\|\bar x_0^\kappa - \mathbf{E}(\bar x_0^\kappa|\mathcal{F}_{]-\infty,0]}^{\bar y,\kappa})\right\|^2
	\right).
$$
The proof is completed by following the same steps as in the proof of 
lemma \ref{lem:manywieners}, and noting that the Wiener process $B$ enters 
linearly in the expression for $\bar x_0^\kappa$.  \qquad
\end{proof}

\begin{corollary}
The filter achieves maximal accuracy if and only if the matrix $H(\lambda 
I-\bar A)^{-1}D$ has linearly independent columns for all 
$\lambda\in\mathbb{C}$ with $\mathrm{Re}\,\lambda>0$.
\end{corollary}

\begin{proof}
Immediate from theorem \ref{thm:KS}. \qquad
\end{proof}

\begin{corollary}
The filter achieves maximal accuracy if and only if the matrix $H(\lambda 
I-A)^{-1}D$ has linearly independent columns for all 
$\lambda\in\mathbb{C}$ with $\mathrm{Re}\,\lambda>0$.
\end{corollary}

\begin{proof}
By \cite[proposition 2]{Gil69}, $H(\lambda I-A)^{-1}D$ has linearly 
independent columns iff $H(\lambda I-A+KH)^{-1}D$ has linearly
independent columns, for any matrix $K$.
\qquad
\end{proof}

\bibliographystyle{siam}
\bibliography{ref}

\end{document}